\documentclass[11pt]{amsart}
\usepackage{amsmath,amssymb,amsfonts,amsthm,amsopn}
\usepackage{graphics}

 \def\Spnr{Sp(d,\R)}
 \def\Gltwonr{GL(2d,\R)}
\newcommand{\fiola}{FIO(\Xi ,s)}

\newcommand{\tfa}{time-frequency analysis}

\newcommand{\stft}{short-time Fourier transform}

\newcommand{\fif}{if and only if}
\newcommand{\tfs}{time-frequency shift}

\newcommand{\modsp}{modulation space}
\newcommand{\psdo}{pseudodifferential operator}

\newtheorem{theorem}{Theorem}[section]
\newtheorem{lemma}[theorem]{Lemma}

\newtheorem{proposition}[theorem]{Proposition}
\newtheorem{definition}[theorem]{Definition}
\newtheorem{cor}[theorem]{Corollary}
\newtheorem{remark}[theorem]{Remark}

\newcommand{\beqa}{\begin{eqnarray*}}
\newcommand{\eeqa}{\end{eqnarray*}}

\newcommand{\field}[1]{\mathbb{#1}}
\newcommand{\bR}{\field{R}}        
\newcommand{\bN}{\field{N}}        
\newcommand{\bZ}{\field{Z}}        
\newcommand{\bC}{\field{C}}        
\def\G{\mathcal{G}}

\def\la{\lambda}

\def\cF{\mathcal{F}}              
\def\cS{\mathcal{S}}

\def\cB{\mathcal{B}}

\def\cG{\mathcal{G}}

\def\cA{\mathcal{A}}
\def\cJ{\mathcal{J}}

 \def\cL{\mathcal{L}}
\def\cC{\mathcal{C}}

\def\cO{\mathcal{O}}

\def\a{\aleph}

\def\rd{\bR^d}

\def\rdd{{\bR^{2d}}}
\def\zdd{{\bZ^{2d}}}

\def\lrd{L^2(\rd)}

\def\zd{\bZ^d}

\def\intrd{\int_{\rd}}

\def\R{\right)}
\def\l{\langle}

\def\<{\left<}
\def\>{\right>}

\def\inv{^{-1}}

\def\mv1{M_v^1}

\def\Mmpq{M_m^{p,q}}
\def\phas{(x,\o )}
\def\mn{(m,n)}
\def\mn'{(m',n')}

\def\Spnr{Sp(d,\R)}

\def\o{\eta}
\def\a{\alpha}

\def\N{\mathbb{N}}
\def\R{\mathbb{R}}
\def\Ren{\mathbb{R}^d}

\def\sch{\mathcal{S}}

\def\Fur{\mathcal{F}}

\def\f{\varphi}

\def\Sn2{S_{2}(L^{2}(\Ren))}
\def\S1{S_{1}(L^{2}(\Ren))}
\def\sig00{\sigma_{0,0}}

\def\la{\langle}
\def\ra{\rangle}

\newcommand{\A}{\mathcal{A}}

\begin{document}
\begin{abstract}
We construct a one-parameter family of algebras $\fiola , s\geq 0,
$ consisting of Fourier integral operators. We derive boundedness
results, composition rules, and the spectral invariance of the
operators in $\fiola $. The operator algebra is defined by the
decay properties of an associated Gabor matrix around the graph of
the canonical transformation.
\end{abstract}

\title{The Wiener Property for a Class of Fourier Integral Operators}

\author{Elena Cordero, Karlheinz Gr\"ochenig, Fabio Nicola,  and Luigi Rodino}
\address{Department of Mathematics,
University of Torino, via Carlo Alberto 10, 10123 Torino, Italy}
\address{Faculty of Mathematics,
University of Vienna, Nordbergstrasse 15, A-1090 Vienna, Austria}
\address{Dipartimento di Matematica,
Politecnico di Torino, corso Duca degli Abruzzi 24, 10129 Torino,
Italy}
\address{Department of Mathematics,
University of Torino, via Carlo Alberto 10, 10123 Torino, Italy}

\email{elena.cordero@unito.it}
\email{karlheinz.groechenig@univie.ac.at}
\email{fabio.nicola@polito.it}
\email{luigi.rodino@unito.it}
\thanks{K.\ G.\ was
  supported in part by the  project P22746-N13  of the
Austrian Science Fund (FWF)}

\maketitle
\section{Introduction}
Wiener's lemma, in its original version \cite{A86}, \cite{A87}, is a classical statement about absolutely convergent series. In a more general setting, Wiener's lemma represents now one of the driving forces in the development of Banach algebra theory.\par
In this paper we will  consider  algebras of Fourier integral
operators (FIOs) and their properties.
 Let us first fix our attention on pseudodifferential operators, which we may express in the Kohn-Nirenberg form
\begin{equation}\label{unoN}
\sigma (x,D)f(x)=\int e^{2\pi i x\eta}
\sigma (x,\eta)\hat{f}(\eta)\,d\eta \, .
\end{equation}
The best known
result about Wiener's property for pseudodifferential operators is
maybe that in \cite{A8}, see also the subsequent contributions of
\cite{bony-chemin,ueberberg}. It concerns symbols $\sigma$ in the
H\"ormander's class $S^{0}_{0,0}(\rdd)$, i.e., smooth functions on
$\rdd$ such that,  for every multi-index $\alpha$ and every
$z\in\rdd$,
\begin{equation}\label{treS}
|\partial^{\alpha}_{z}\sigma(z)|\leq C_{\alpha}.
\end{equation}
The corresponding \psdo s   form a subalgebra of
$\mathcal{L}(L^{2}(\rd))$, usually denoted by $L^{0}_{0,0}$. The
standard symbolic calculus concerning the  principal part of symbols
of products does not hold, nevertheless Wiener's lemma is
 still valid. Namely, if $\sigma (x,D)
$ is invertible in $\mathcal{L}(L^{2}(\rd))$, then its inverse is
again a pseudodifferential operator with symbol in
$S^{0}_{0,0}(\rdd)$, hence  belonging to $L^{0}_{0,0}$. In the
absence of a symbolic calculus, such a version of Wiener's lemma
seems to be the minimal property required of any reasonable
algebra of pseudodifferential operators. A subalgebra $\cA $ of
operators in $\cL (\lrd )$ that satisfies Wiener's lemma and is
thus  closed under inversion, is usually called \emph{spectrally
invariant} or inverse-closed or sometimes also a Wiener algebra.
See~\cite{A-G} for a survey of the theory of spectral invariance.

 From the point of
view of time-frequency analysis and signal processing, which we are going to adopt in
this paper, Wiener's lemma provides an  important justification of the
engineering practice  to model  $\sigma(x,D)^{-1}$
as an almost diagonal matrix (this  is a peculiar property of
pseudodifferential operators, see Theorem \ref{CR1} in the
sequel). \par
Actually, in the applications to signal processing, the symbol
$\sigma(x,\eta)$ is not always  smooth,  and it
is convenient to use  some generalized version of
$S^{0}_{0,0}(\rdd)$~\cite{str06}. Let us recall some  results in this connection. To
give a unified presentation, we use the modulation spaces $M^{p,q}_m$ introduced
by Feichtinger, cf.~\cite{F1,book}. 
See Section 2.1 for the definition.  We are
particularly interested in the so-called Sj\"ostrand class
\begin{equation}\label{quattro}
S_w=M^{\infty,1}(\rdd)
\end{equation}
and the related scale of spaces
\begin{equation}\label{cinque}
S^s_w=M^{\infty,\infty}_{1\otimes v_s}(\rdd),\quad v_s(z)=\langle
z\rangle ^s=(1+|z|^2)^{s/2}, \,\,z\in\rdd,
\end{equation}
with  the  parameter $s\in [0,\infty )$. Defining
$S^\infty_w=\bigcap_{s\geq 0} S^s_w$, we recover the H\"ormander
class $S^\infty_w=S^0_{0,0}(\rdd)$, whereas for $s\to 2d$ the
symbols in $S^s_w$ have a smaller  regularity, until in  the
maximal space $S_w $ even  differentiability is lost.

\begin{theorem}[\cite{A76,GR}]\label{CR2}
The  pseudodifferential operators with a  symbol in $S_w$ form a
Wiener subalgebra of $\mathcal{L}(L^2(\rd))$, the so-called
Sj\"ostrand  algebra. The symbol classes  $S^s_w$ with  $s>2d$  provide a
scale of Wiener subalgebras of the Sj\"ostrand  algebra, and their
intersection coincides with $L^0_{0,0}$.
\end{theorem}\par\bigskip
In this paper we construct and  investigate  Wiener subalgebras
consisting of Fourier integral operators  and generalize
Sj\"ostrand's theory in ~\cite{A76,Sjo95} to FIOs. We will  consider FIOs of type I, that is
\begin{equation}\label{sei}
T f(x)=T_{I,\Phi,\sigma}f(x)=\int_{\rd} e^{2\pi i
  \Phi(x,\eta)}\sigma(x,\eta)\widehat{f}(\eta)\,d\eta \, ,
\end{equation}
where we first assume $\sigma\in S^0_{0,0}(\rdd)$. For the
real-valued phase $\Phi$ we assume that $\partial^\alpha\Phi(z)\in
S^0_{0,0}(\rdd)$ for $|\alpha|\geq2$ and that a standard
non-degeneracy  condition is satisfied,
cf. Section 2. When $\Phi(x,\eta)=x\eta$ we recapture the
pseudodifferential operators in the Kohn-Nirenberg form. The
$L^2$-adjoint of a FIO of type $I$ is a FIO of type II
\begin{equation}\label{sette}
Tf(x)=T_{II,\Phi,\tau}f(x)=\int_{\rdd} e^{-2\pi i[\Phi(y,\eta)-x\eta]}\tau(y,\eta) f(y)\, dy\, d\eta.
\end{equation}
 For the $L^2$-boundedness of such FIOs of types I and II see for
example \cite{AS78}.

It is worth to observe that in contrast   to the standard setting
of H\"or\-man\-der~\cite{B36},  we argue globally on $\rd$,  our
basic examples being the propagators for Schr\"odinger-type
equations. The second remark is that operators of the reduced form
\eqref{sei} or \eqref{sette} do not form an algebra,  quite in
line with  the calculus of \cite{B36}. The composition of FIOs
requires heavier machinery and is  addressed,  for example,  in
\cite{B33} in  the case when symbol and phase belong to the more
restrictive Shubin class of \cite{B43}.\par As minimal objective
of the present paper we want to present a cheap definition of a
subalgebra of $\mathcal{L}(L^2(\rd))$ containing FIOs of type I,
II, and hence $L^0_{0,0}$, and prove the  Wiener property for this
class.\par As a  more  ambitious objective we will  extend our
analysis to the case when $\sigma$ in \eqref{sei} belongs to the
symbol class $S^s_w$, and we will  define  a corresponding scale
of Wiener algebras of FIOs (we will not treat the full Sj\"ostrand
algebra  in this paper).  The new  algebras of FIOs will be
constructed by means of Gabor frames and  the decay properties of
the corresponding Gabor
matrix 
 outside the graph of a symplectic map $\chi$. For FIOs of type I such
a decay 
was already pointed out in
\cite{CNG,fio1} with applications to boundedness properties and
numerical analysis in \cite{fio3}.
\par\bigskip
For the formulation of the results we now  introduce  the basic
notions of  time-frequency analysis  and
refer  to Section 2 for  details.
The most suitable  representation  for our purpose  is the  \stft,
where the localization on the time-frequency plane $\rdd$ occurs on
the unit scale both in time and in frequency.
 For a point $z=\phas\in\rdd$ and a function $f$ on $\rd$,
we denote the  time-frequency shifts (or phase-space
shifts) by $$\pi(z)=M_\eta T_x f(t)= e^{2\pi i t\eta} f(t-x),\quad
\mbox{where}\quad t\eta=t\cdot\eta=\sum_{i=1}^d t_i\eta_i \, .$$

The \stft\, (STFT) of a function/distribution $f$ on $\rd$ with respect to a Schwartz window function $g\in\cS(\rd)\setminus\{0\}$ is defined by

\begin{equation}\label{stft}
V_g f(x,\o)=\langle f,M_{\o}T_x g\rangle=\langle f,\pi(z)g
\rangle= \int _{\rd} f(t)\overline{g(t-x)}e^{-2\pi it\o}\,dt,
\end{equation}
for $z=(x,\eta) \in\rdd$.
We can now  define the generalized Sj\"ostrand class  $S^s_w$ in
\eqref{cinque}  as the space of distributions $\sigma\in\cS'(\rdd)$
such that
\begin{equation}
 |\langle \sigma ,\pi(z,\zeta )g\rangle|\leq C\langle \zeta
 \rangle^{-s},\quad \forall z,\zeta \in \rdd 
 \end{equation}
 for some constant $C>0$, whereas  $\sigma$ is in the  Sj\"ostrand class $S_w$  if   $$\int_{\rd}\sup_{z\in\rd}|\langle \sigma ,\pi(z,\zeta )g\rangle|
 \,d\zeta<\infty.$$

 \par

For the discrete description of function spaces and operators we use
Gabor frames. Let $\Lambda=A\zdd$  with $A\in GL(2d,\R)$ be a lattice
of the time-frequency plane.
 The set  of
time-frequency shifts $\G(g,\Lambda)=\{\pi(\lambda)g:\
\lambda\in\Lambda\}$ for a  non-zero $g\in L^2(\rd)$ is called a
Gabor system. The set $\G(g,\Lambda)$   is
a Gabor frame, if there exist
constants $A,B>0$ such that
\begin{equation}\label{gaborframe}
A\|f\|_2^2\leq\sum_{\lambda\in\Lambda}|\langle f,\pi(\lambda)g\rangle|^2\leq B\|f\|^2_2\qquad \forall f\in L^2(\rd).
\end{equation}
 Gabor frames allow us to discretize  any continuous operator from
 $\cS(\rd)$ to $\cS'(\rd)$ into  an infinite matrix that captures the  properties of the original operator. \par

For  the case of
pseudodifferential operators  the Gabor discretization provides an
equivalent characterization of the Sj\"ostrand algebra in Theorem \ref{CR2}.
\begin{theorem}[\cite{charly06,GR}]\label{CR1}  Assume  that $\G(g,\Lambda)$ is a
 frame for $L^2(\rd)$ with  $g\in\cS(\rd)$ and fix $s>2d$. Then the
 following statements  are equivalent
for a distribution $\sigma\in\cS'(\rd)$:
\par {\rm
(i)} $\sigma\in S^s_w$.\par
{\rm (ii)} There exists $C>0$ such that
\begin{equation}\label{unobis2s} |\langle \sigma (x,D) \pi(\lambda)
g,\pi(\mu)g\rangle|\leq C\langle\mu-\lambda\rangle^{-s},\qquad \forall
\lambda,\mu\in \Lambda.
\end{equation}
Hence, the assumption $\sigma\in S^\infty_w=S^0_{0,0}(\rd)$ is
equivalent \eqref{unobis2s} being  satisfied for all $s\geq 0$.\par
Moreover $\sigma\in S_w$ \fif\ there exists a sequence $h\in \ell
^1(\Lambda )$, such that  $|\langle \sigma (x,D) \pi(\lambda)
g,\pi(\mu)g\rangle|\leq h(\lambda-\mu)$.
\end{theorem}

We refer to the matrix of an operator $T$ with respect to a Gabor frame as
the \emph{Gabor matrix} of $T$.
The above  theorem gives a precise meaning to the statement that the Gabor
matrix of a  pseudodifferential operators  is  almost diagonal, or
that \psdo s are almost diagonalized by Gabor frames. \par We
 now  describe our results about FIOs in Section 3.

Roughly speaking,  Fourier integral operators can be defined as follows (Definition \ref{defFIO}).
Consider a bi-Lipschitz canonical transformation  $\chi: \rdd\to\rdd$
(see Definition \ref{de}) and $s > 2d$. Let
$\G(g,\Lambda)$ be a Gabor frame for $L^2(\rd)$ with $g\in \cS (\rd )$. We say that  a
continuous linear operator $T:\cS(\rd)\to\cS'(\rd)$ is in the
class $FIO(\chi,s)$ if its Gabor matrix  satisfies the decay
condition   \begin{equation}\label{asterisco}
|\langle T \pi(\lambda) g,\pi(\mu)g\rangle|\leq {C}\langle \mu-\chi(\lambda)\rangle^{-s},\qquad \forall \lambda,\mu\in \Lambda.
\end{equation}

If $\chi=\mathrm{Id}$, the identity operator, then the corresponding Fourier integral operators are simply pseudodifferential operators.

The  decomposition of a FIO with respect to a Gabor frame  provides a
technique to settle the following issues.
\begin{itemize}
    \item [(i)] \emph{Boundedness of $T$ on $\lrd$} (Theorem \ref{P26}):\\
    If $s>2d$ and $T\in FIO(\chi,s)$, then $T$ can be extended to a bounded operator on $\lrd$.
    \item[(ii)] \emph{The algebra property} (Theorem \ref{prod1}): For $ i=1,2$,  $s>2d,$
$$T^{(i)}\in FIO(\chi_i,s_i)\quad \Rightarrow \quad T^{(1)}T^{(2)}\in
FIO(\chi_1\circ \chi_2, s) \, .$$
 \item[(iii)] \emph{Wiener property} (Theorem \ref{Wineroper}): If $s>2d$, $T\in FIO(\chi,s)$  and  $T$ is invertible on $L^2(\rd)$, then $T^{-1} \in FIO(\chi^{-1},s)$.
\end{itemize}
These three properties can be summarized neatly by saying that the
union $\bigcup _{\chi } FIO(\chi ,s)$ is a Wiener subalgebra of $\cL
(\lrd )$ consisting of FIOs.
 \par
 \medskip\noindent
In Section 4 we return to concrete FIOs  of type I and
II. Denoting by $\chi$ the symplectic transformation related to
a phase $\Phi$, we prove the expected extension of Theorem \ref{CR2} to
FIOs. Namely: A FIO $T$ of type I  as in \eqref{sei} belongs to $FIO(\chi,s)$ for some  $s>2d$,
 if and only if its symbol
$\sigma$ belongs to $S^s_w$ (Theorem \ref{caraI}). We
further prove  that the inverse in $\mathcal{L}(L^2(\rd))$ of an operator of type I is an operator of type II, with symbol belonging to the same class $S^s_w$ (Theorem \ref{Winerfio}). As an example, in Section 5 we treat a Wiener algebra of generalized metaplectic operators.\par\bigskip
Although it is impossible to do justice to the vast literature on
Fourier integral operators, let us mention some of the contributions
that are most related to our
ideas.  From
the formal point of view, our approach is very similar to that in
\cite{bony1,bony3} and \cite{tataru}, where $FIO(\chi,\infty)$ was
treated. Instead of Gabor frames, in \cite{bony1,bony3} partitions of
unity of the Weyl-H\"ormander calculus are used, whereas in
\cite{tataru} the Bargmann transform is the main tool. The boundedness
and composition of FIOs are treated in  \cite{Bou97,
  Bou97french,concetti-garello-toft,concetti-toft,B18,B33}.\par
The  time-frequency analysis of \psdo s was propagated in
\cite{A34,A40,GR,rochberg,A84}. Many aspects of Wiener's lemma and
spectral invariance of operators are surveyed in ~\cite{A-G}.
 \vspace{5mm}

\textbf{Notation.} We write $xy=x\cdot y$ for  the scalar product on
$\Ren$ and $|t|^2=t\cdot t$ for $t,x,y \in\Ren$.

The Schwartz class is denoted by
$\sch(\Ren)$, the space of tempered
distributions by  $\sch'(\Ren)$.   We
use the brackets  $\la f,g\ra$ to
denote the extension to $\sch '
(\Ren)\times\sch (\Ren)$ of the inner
product $\la f,g\ra=\int f(t){\overline
{g(t)}}dt$ on $L^2(\Ren)$. The Fourier
transform is normalized to be ${\hat
  {f}}(\o)=\Fur f(\o)=\int
f(t)e^{-2\pi i t\o}dt$.

For $1\leq p\leq\infty$ and a  weight $m$, the space $\ell^{p}_m
(\Lambda )$
 is  the Banach
space of sequences $a=\{{a}_{\lambda}\}_{\lambda \in \Lambda }$
on a  lattice $\Lambda$, such that
$$\|a\|_{\ell^{p}_m}:=\left(\sum_{\lambda\in\Lambda}
|a_{\lambda}|^p m(\lambda)^p\right)^{1/p}<\infty
$$
(with obvious changes when $p=\infty$).

Throughout the paper, we
shall use the notation
$A\lesssim B$ to express the inequality
$A\leq c B$ for a suitable
constant $c>0$, and  $A
\asymp B$  for the equivalence  $c^{-1}B\leq
A\leq c B$.

 \section{Preliminaries}

 \subsection{Phase functions and canonical transformations}

 \begin{definition}\label{de}
A real phase function $\Phi$ on $\rdd$ is called \emph{tame}, if the
following three properties are satisfied:\par\medskip \noindent
A1. $\Phi\in \cC^{\infty}(\rdd)$;\\
A2. For $z=\phas$,
\begin{equation}\label{phasedecay}
|\partial_z^\a \Phi(z)|\leq C_\a,\quad |\a|\geq 2;\end{equation}
A3. There exists $\delta>0$ such that
\begin{equation}\label{detcond}
   |\det\,\partial^2_{x,\eta} \Phi(x,\o)|\geq \delta.
\end{equation}
\end{definition}
\par
If we set
\begin{equation}\label{cantra} \left\{
                \begin{array}{l}
                y=\nabla_{\eta}\Phi(x,\eta)
                \\
               \xi=\nabla_{x}\Phi(x,\eta), \rule{0mm}{0.55cm}
                \end{array}
                \right.
\end{equation}
we can solve with respect to $(x,\xi)$ by the global inverse
function theorem (see e.g. \cite{krantz}) and  obtain a mapping
$\chi$ defined by $(x,\xi)=\chi(y,\o)$. The canonical transformation
$\chi $ enjoys the following properties:

\par\medskip
\noindent {\it  B1.} $\chi:\rdd\to\rdd$ is smooth, invertible,  and
preserves the symplectic form in $\rdd$, i.e., $dx\wedge d\xi= d
y\wedge d\eta$; $\chi $ is a \emph{symplectomorphism}.
 \\
{\it  B2.} For $z=(y,\eta)$,
\begin{equation}\label{chistima}
|\partial_z^\a \chi(z)|\leq C_\a,\quad |\a|\geq 1;\end{equation}
{\it B3}. There exists $\delta>0$ such that, for
$(x,\xi)=\chi(y,\eta)$,
\begin{equation}\label{detcond2}
   |\det\,\frac{\partial x}{\partial y}(y,\eta)|\geq \delta.
\end{equation}
Conversely,  to every transformation $\chi$ satisfying {\it B1},
{\it B2}, {\it B3} corresponds a tame phase $\Phi$, uniquely
determined up to a constant. This can be easily proved by
\eqref{detcond2}, the global inverse function theorem
\cite{krantz} and using the pattern of \cite[Theorem
4.3.2.]{Rodino} (written for the local case). \par
From now on we shall
define by $\Phi_\chi$ the phase function (up to constants)
corresponding to the canonical transformation $\chi$.

Observe that {\it B1} and {\it B2} imply that $\chi$ and
$\chi^{-1}$ are globally Lipschitz.
This property implies that
$$
\langle w-\chi (z) \rangle \asymp \langle \chi \inv (w) - z \rangle
\qquad w,z\in \rdd \, ,
$$
which we will use frequently.
  Moreover, if $\chi$ and
$\tilde{\chi}$ are two transformation satisfying {\it B1} and {\it
B2}, the same is true for $\chi\circ \tilde{\chi}$, whereas the
additional property {\it B3} is  not necessarily  preserved, even if $\chi$ and
$\tilde{\chi}$ are linear. This reflects  the lack of the algebra
property of  the corresponding FIOs of type I; see Section
\ref{metapoper} below.

\subsection{Time-frequency concepts}\label{tfc}
We recall the basic
concepts  of \tfa\ and  refer the  reader to \cite{book} for the full
details. 
Consider a distribution $f\in\cS '(\rd)$
and a Schwartz function $g\in\cS(\rd)\setminus\{0\}$ (the so-called
{\it window}).
The short-time Fourier transform of $f$ with respect to $g$ was
defined in  \eqref{stft} by $V_gf (z) = \langle f, \pi (z)g\rangle
$. The  \stft\ is well-defined whenever  the bracket $\langle \cdot , \cdot \rangle$ makes sense for
dual pairs of function or distribution spaces, in particular for $f\in
\cS ' (\rd )$ and $g\in \cS (\rd )$ or $f,g\in\lrd$.
We recall the covariance formula for the \stft\  that will be used in
the sequel (Section \ref{metapoper}):
  \begin{equation}
    \label{eql2}
   V_{ {g}}(M_\xi T_y{f})(x,\o)  = e^{-2\pi
     i(\o-\xi)y}(V_gf)(x-y,\o-\xi), \qquad x, y,\omega , \xi \in \Ren
   .
  \end{equation}

\par
The symbol spaces are provided by the  modulation spaces. These were introduced by Feichtinger in the 80's (see the original paper \cite{F1}) and now are well-known in the framework of time-frequency analysis.
For the fine-tuning of decay properties in the definition of \modsp s
we use  weight functions of polynomial growth.  For $s\geq 0$
we set $v\phas=v_s\phas=\la
\phas\ra^s=(1+|x|^2+|\o|^2)^{s/2}$ and  denote by
$\mathcal{M}_v(\rdd)$ the space of $v$-moderate weights on $\rdd$;
these  are measurable functions $m>0$ satisfying $m(z+\zeta)\leq C
v(z)m(\zeta)$ for every $z,\zeta\in\rd$.
In particular, $v_s (z) \inv = \langle z\rangle ^{-s}$ is
$v_s$-moderate. The corresponding inequality $\langle z+w\rangle ^{-s} \leq \langle
z\rangle ^{-s} \langle w\rangle ^s$ is also called Peetre's
inequality.
\par

Let $g$ be a non-zero Schwartz  function. For $1\leq p,q \leq \infty$ and $m\in
\mathcal{M}_v(\rdd)$ the modulation space $M ^{p,q}_m(\R^d)$ is the space of
distributions $f\in\cS'(\rd)$ such that their STFTs belong to the
space $L^{p,q}_m(\rdd ) $ with  norm
\[
\|f\|_{M ^{p,q}_m(\R^d)}:=\|V_g f\|_{L^{p,q}_m(\R^{2d})}=\left(\intrd\left(\intrd|V_g f\phas|^p m\phas^p dx\right)^{\frac q p}d\eta\right)^{\frac1q}.
\]
This definition  does not depend on the choice of
the window $g\in \cS (\rd ), g \neq 0$, and different windows yield
equivalent norms on $\Mmpq$~\cite[Thm.~11.3.7]{book}. Moreover, the space of admissible windows can be enlarged to $M^1_v(\rd )$.
The symbol spaces we shall  be mainly concerned with are
$S^s_w=M^{\infty,\infty}_{1\otimes v_s}(\rdd)$ with the  norm
$$\|\sigma\|_{S^s_w}=\sup_{z\in\rdd}\sup_{\zeta\in\rdd} |V_\Psi \sigma
(z,\zeta)|\,\langle \zeta \rangle ^s,
$$
where $\Psi\in\cS(\rdd)\setminus\{0\}$.
\par
The H\"ormander  symbol class $S^0_{0,0}(\rdd)$  can be characterized
by means of modulation spaces as
follows, see for example ~\cite{GR}: $$S^0_{0,0}=\bigcap_{s\geq 0}
S^s_w.$$

\subsubsection{Gabor frames}
Fix a function $g\in\lrd$ and a lattice
$\Lambda =A\zd$, for
$A\in GL(2d,\R)$.
The Gabor system $\cG(g,\Lambda)=\{\pi(\lambda) g:\,\lambda\in\Lambda\}$ is a Gabor frame
if there exist constants $A,B>0$ such that \eqref{gaborframe} is satisfied.
We define the coefficient operator $C_g$, which
maps functions to sequences as follows:
\begin{equation}\label{analop}
    (C_gf)_{\lambda}:=\la
    f,\pi(\lambda)g\ra,\quad \lambda\in \Lambda,
\end{equation}
the synthesis operator
\begin{equation}
  \label{eq:c6}
D_g{c}:=\sum_{\lambda\in \Lambda}
c_{\lambda} \pi(\lambda)g,\quad
c=\{c_{\lambda}\}_{\lambda\in
\Lambda}
\end{equation}
and the Gabor frame operator
\begin{equation}\label{Gaborop}
    S_g f:=D_g C_g f=\sum_{\lambda\in \Lambda}\la f,\pi(\lambda)g\ra \pi(\lambda)g.
\end{equation}

Equivalently,  the set $\cG(g,\Lambda)$  is called a
Gabor frame for the Hilbert space
$\lrd$,  if $S_g$ is a bounded and
 invertible operator on $\lrd$. If $\cG(g,\Lambda)$ is a Gabor frame for $\lrd$,
  then the so-called \emph{dual window}
   $\gamma=S_g^{-1} g$ is well-defined and the
   set $\cG(\gamma,\Lambda)$  is a frame (the so-called
   canonical dual frame of $\cG(g,\Lambda)$). Every $f\in \lrd$ possesses the frame expansion
\begin{equation}\label{frame}
    f= \sum_{\lambda\in \Lambda}\la f,\pi(\lambda)g\ra \pi(\lambda)\gamma= \sum_{\lambda\in \Lambda}\la f,\pi(\lambda)\gamma\ra \pi(\lambda)g
\end{equation}
with unconditional convergence in
$\lrd$, and norm equivalence
$$\|f\|_{L^2}\asymp \|C_g f\|_{\ell^2}\asymp \|C_\gamma f\|_{\ell^2}.
$$
Gabor frames give the following characterization of the Schwartz space
$\cS(\rd)$  and of the modulation spaces $M^{p}_m(\rd)$, $1\leq
p\leq\infty$. If  $g\in \cS (\rd )$ and $\cG (g,\Lambda )$ is a frame, then
\begin{align}
 \label{osservazione0}
f\in\cS(\rd) & \Leftrightarrow \sup_{\lambda\in\Lambda}\la \lambda\ra^N
|\langle f,\pi(\lambda)g\rangle|<\infty \quad\forall N\in\bN,\\
\label{osservazione0bis}
f\in M^{p}_m(\rd) &\Leftrightarrow \Big(\sum_{\lambda\in\Lambda}|\langle f,\pi(\lambda)g\rangle|^p m(\lambda)^p\Big)^{1/p}<\infty,
\end{align}
These results are contained in
\cite[Ch.~13]{book}. In
particular, if $\gamma=g$, then the frame is called
\emph{Parseval}  frame and
the expansion \eqref{frame} reduces to
\begin{equation}\label{parsevalframe}
    f= \sum_{\lambda\in \Lambda}\la f,\pi(\lambda)g\ra \pi(\lambda)g.
\end{equation}
 We may take the existence of Parseval frames with $g\in \cS (\rd )$
 for granted.   From now on we  work with Parseval frames,   so that
 we will not have to deal with the dual window $\gamma$. Let us
 underline that the properties of FIOs written for  Parseval  frames
 work exactly the same with general Gabor frames with dual windows
 $\gamma$ different from $g$.


\section{A Wiener Algebra of Fourier Integral Operators}
We first present an equivalence between continuous decay conditions
and the decay of the discrete Gabor
matrix for a linear operator $\cS(\rd)\to\cS'(\rd)$.

\begin{theorem}\label{cara}
Let $T$ be a continuous linear operator $\cS(\rd)\to\cS'(\rd)$ and
$\chi$ a canonical transformation which satisfies {\it B1} and {\it
B2} of Definition \ref{de}. Let $\G(g,\Lambda)$ be a Parseval frame with
$g\in\cS(\rd)$ and $s\geq 0$. Then the following properties are
equivalent. \par {\rm
(i)} There exists $C>0$ such that
\begin{equation}\label{unobis}
|\langle T \pi(z) g,\pi(w)g\rangle|\leq {C}\langle w-\chi(z)\rangle^{-s},\qquad \forall z,w\in \rdd.
\end{equation}
\indent{\rm (ii)}
There exists $C>0$ such that
\begin{equation}\label{unobis2}
|\langle T \pi(\lambda) g,\pi(\mu)g\rangle|\leq {C}\langle \mu-\chi(\lambda)\rangle^{-s},\qquad \forall \lambda,\mu\in \Lambda.
\end{equation}
\end{theorem}
\begin{proof}  The implication  ${\rm (i)\Longrightarrow{\rm  (ii) }}$
  is obvious.\par
${\rm (ii)}\Longrightarrow{\rm  (i) }$. The argument  is borrowed from
the proof of a similar result for pseudodifferential operators
in \cite[Theorem 3.2]{charly06}.  Let $C$ be a relatively compact
fundamental domain of the lattice $\Lambda$. Given $z,w\in\rdd$,
we can write $w=\lambda+u$, $z=\mu+u'$ for unique
$\lambda,\mu\in\Lambda$, $u,u'\in C$, and
\[
|\langle T\pi(z)g,\pi(w) g \rangle|=|\langle
T\pi(\mu)\pi(u')g,\pi(\lambda) \pi(u)g \rangle|.
\]
Now we expand  $\pi(u)g=\sum_{\nu\in\Lambda} \langle \pi(u)
g,\pi(\nu) g\rangle \pi(\nu) g$, and likewise $\pi(u')g$. Since $V_gg
\in \cS (\rdd )$, the coefficients in this expansion satisfy
\[
\sup_{u\in C}|\langle \pi(u)g,\pi(\nu)g\rangle|= \sup_{u\in C}|V_g
g(\nu-u)|\lesssim \sup_{u\in C}\langle \nu-u\rangle ^{-N}\lesssim
(\sup_{u\in C}\langle u\rangle ^{N})\langle \nu\rangle
^{-N}\lesssim \langle \nu\rangle ^{-N}
\]
 for every $N$.
Using  \eqref{unobis2}, we now  obtain that
\begin{align*}
|\langle T\pi(z)g,\pi(w) g \rangle|&\leq\sum_{\nu,\nu'\in\Lambda}|\langle T \pi(\mu+\nu') g,\pi(\lambda+\nu)g\rangle|
|\langle \pi(u')g,\pi(\nu')g\rangle| |\langle \pi(u)g,\pi(\nu)g\rangle| \\
&\lesssim  \sum_{\nu,\nu'\in\Lambda}\langle \lambda+\nu-\chi(\mu+\nu')\rangle^{-s}\langle \nu'\rangle ^{-N} \langle \nu\rangle^{-N}.
\end{align*}
Since $v_s(z)\inv = \langle z\rangle ^{-s}$ is $v_s$-moderate, we
majorize the main  term of the sum as
\begin{align*}
  \langle \lambda+\nu-\chi(\mu+\nu')\rangle^{-s} & =   \langle \lambda
  - \chi (\mu ) +\nu-\chi(\mu+\nu') + \chi (\mu ) \rangle^{-s}\\
& \leq \langle \lambda
-\chi(\mu)\ra^{-s}\langle
\chi(\mu+\nu')-\chi(\mu)-\nu\rangle^s\lesssim \langle \lambda
-\chi(\mu)\ra^{-s}\la \nu'\ra^s\la\nu\ra^s \, ,
\end{align*}
since $\chi (\mu +\nu ') - \chi (\mu ) = \cO (\nu ')$ by the Lipschitz
property of $\chi $.

Hence,
\begin{equation*}
|\langle T\pi(z)g,\pi(w) g \rangle|\lesssim \langle \lambda
-\chi(\mu)\ra^{-s} \sum_{\nu,\nu'\in\Lambda}\langle
\nu'\rangle^{s-N}\langle \nu\rangle^{s-N}\lesssim \langle \lambda
-\chi(\mu)\ra^{-s},
\end{equation*}
for $N\in\N$ large enough.

Finally, with   $\lambda=w-u$, $\mu=z-u'$, we apply the above estimate
again   and obtain
\begin{align*}
\langle \lambda - \chi (\mu )\rangle ^{-s} &= \langle w-u  - \chi
(z-u'  ) -\chi (z ) + \chi (z) \rangle ^{-s}  \\
&\leq \langle w - \chi (z )\rangle ^{-s} \langle u+ \chi (z-u') -\chi
(z) \rangle ^s \\
&\lesssim \langle w - \chi (z )\rangle ^{-s} \langle u\rangle
^s\langle u '\rangle ^s \lesssim  \langle w - \chi (z
)\rangle ^{-s} \, ,
\end{align*}
since $\sup _{u\in C} \langle u \rangle ^s < \infty$.  Thus we have
proved that $ |\langle T\pi(z)g,\pi(w) g \rangle| \lesssim \langle w -
\chi (z )\rangle ^{-s}$.
\end{proof}

Inspired by  the characterization of Theorem \ref{cara}, we now
 define a   class of FIOs associated to a canonical
transformation  $\chi$. In view of the
equivalence \eqref{unobis} and \eqref{unobis2}, we focus on the
 decay of the \emph{discrete} Gabor matrix.

\begin{definition}\label{defFIO}
Let $\chi$ be a  transformation satisfying {\it B1} and
{\it B2}, and $s\geq0$. Fix $g\in\cS(\rd)\setminus\{0\}$ and let
$\G(g,\Lambda)$ be a Parseval frame for $L^2(\rd)$. We say that  a
continuous linear operator $T:\cS(\rd)\to\cS'(\rd)$ is in the
class $FIO(\chi,s)$, if its Gabor matrix  satisfies the decay
condition
\begin{equation}
  \label{unobis2a}
|\langle T \pi(\lambda) g,\pi(\mu)g\rangle|\leq {C}\langle \mu-\chi(\lambda)\rangle^{-s},\qquad \forall \lambda,\mu\in \Lambda.
\end{equation}
The class $\fiola = \bigcup _{\chi } FIO (\chi ,s)$ is  the union of
these classes where  $\chi$ runs over the  set of all transformations satisfying $B1,B2$.
\end{definition}
Note that we do not require the
assumption {\it B3}.

We first observe that this definition  does not depend on the choice
of the Gabor  frame.

\begin{lemma}
  The definition of $FIO(\chi ,s)$ is independent of the Gabor frame
  $\cG (g, \Lambda ) $.
\end{lemma}

\begin{proof}
 Let $\G(\f,\Lambda ' )$ be a Gabor frame with a window
 $\f\in\cS(\rd) $ and a possibly different lattice $\Lambda '$. 
As in the proof of Theorem~\ref{cara} we expand
 $\pi(\lambda)\f =\sum_{\nu\in\Lambda} \la\pi(\lambda)\f,\pi(\nu)g\ra
 \pi(\nu)g$ with convergence in $\cS(\rd)$ and likewise $\pi (\mu )\f
 $, where $\lambda ,\mu \in \Lambda '$. Consequently
$T\pi(\lambda)\f =\sum_{\nu\in\Lambda} \la\pi(\lambda)\f,\pi(\nu)g\ra
T\pi(\nu)g$  converges weak$^*$ in $\cS'(\rd)$
and the following identity  is well-defined:
\begin{align*}
\langle T\pi(\lambda)\f,\pi(\mu)\f\ra &=\sum_{\nu\in\Lambda} \la\pi(\lambda)\f,\pi(\nu)g\ra \la T\pi(\nu)g,\pi(\mu)\f\ra\\
&=\sum_{\nu\in\Lambda} \sum_{\nu'\in\Lambda}\la\pi(\lambda)\f,\pi(\nu)g\ra \la T\pi(\nu)g,\pi(\nu')g\ra \la
\pi(\mu)\f,\pi(\nu')g\ra
\end{align*}
Since $\f , g \in \cS (\rd )$, the characterization ~\eqref{osservazione0} and the
covariance property~\eqref{eql2} imply that
$|\la\pi(\lambda)\f,\pi(\nu)g\ra|=|\la\f,\pi(\nu-\lambda)g\ra|\lesssim
\la \nu-\lambda\ra^{-N}$ for $\nu \in \Lambda , \lambda \in \Lambda '$
and every $N\geq 0$. 
After substituting these estimates and choosing $N$ large enough, we obtain the majorization
\begin{align*}
|\langle T\pi(\lambda)\f,\pi(\mu)\f\ra |&\lesssim \sum_{\nu\in\Lambda} \sum_{\nu'\in\Lambda} \la \nu-\lambda\ra^{-N}
\la \nu'-\chi(\nu)\ra^{-s}\la \nu'-\mu\ra^{-N }\\
&\lesssim \sum_{\nu\in\Lambda} \la \nu-\lambda\ra^{-N}
\la \mu-\chi(\nu)\ra^{-s}\asymp \sum_{\nu\in\Lambda} \la \nu-\lambda\ra^{-N}
\la \chi^{-1}(\mu)-\nu\ra^{-s}\\
&\lesssim \la \chi^{-1}(\mu)-\lambda\ra^{-s}\asymp \la
\mu-\chi(\lambda)\ra^{-s}\, .
\end{align*}
 \end{proof}

As in \cite{GR}, the definition of classes of  operators by
their  Gabor matrices facilitates the investigation of their basic
properties. In line with Sj\"ostrand's original program  we next  derive the boundedness, the
composition rules, and properties of the inverse operator in the
classes $FIO(\chi ,s)$.

For  $s>2d$  the class $FIO(\chi,s)$ possesses many desired
properties. 
\begin{theorem}\label{P26}
Let $s>2d$ and  $T\in FIO(\chi,s)$. Then $T$ extends to a bounded
operator on $M^p(\rd )$, $1\leq p\leq\infty$, and in particular on
$L^2(\rd) = M^2(\rd )$. 
\end{theorem}
\begin{proof}
Let $\G(g,\Lambda)$ be a Parseval frame with
$g\in\cS(\rd)$. Since the frame operator $S_g = D_g C_g$ is the
identity operator, we can write $T$ as  $T= D_g C_g T D_g
C_g$, where $D_g$ and $C_g$ are the synthesis and coefficient
operators of \eqref{analop} and \eqref{eq:c6}.
Since $\cG (g,\Lambda )$ is a frame and $g\in \cS (\rd )$, $C_g$ is
bounded from $M^p(\rd )$ to $\ell ^p(\Lambda )$ and $D_g = C_g^*$ is
bounded from  $\ell ^p(\Lambda )$  to $M^p(\rd )$. (For $p=2$ this is
contained in the definition of a frame, for $p\neq 2$ this is slightly
less obvious and stated in~\cite[Ch.~12.2]{book}.)

The operator $C_g T D_g$ maps sequences to sequences, and its matrix
$K$ is precisely the Gabor matrix of $T$, namely,  $K_{\mu, \lambda}=\la
T\pi(\lambda)g,\pi(\mu)g\ra$.
Since by assumption \linebreak  $|\la
T\pi(\lambda)g,\pi(\mu)g\ra |\lesssim \langle \mu - \chi (\lambda
)\rangle ^{-s}$ and $s>2d$, Schur's test implies that the  matrix $K$
representing $C_g T D_g$ is bounded on $\ell ^p(\Lambda
)$. Consequently $T$ is bounded on $M^p(\rd )$ for $1\leq p \leq
\infty $.
\end{proof}

\begin{remark} As it is clear from the preceeding proof, the boundedness of the operator $T\in FIO(\chi,s)$, $s>2d$, fails
in general on the modulation spaces $M^{p,q}(\rdd)$, with $p\not=q$. Indeed, in this case the change of variables $u=\chi(z)$ is not allowed. A concrete counter-example is the FIO of type I in \eqref{moltiplicatore} below.
\end{remark}

Next, we show  that the class $\fiola $ for   $s>2d$ is an  algebra.

\begin{theorem}\label{prod1}
If $T^{(i)}\in FIO(\chi_i,s_i)$ with $s_i>2d$, $i=1,2$; then the
composition  $T^{(1)}T^{(2)}$ is in $FIO(\chi_1\circ \chi_2, s)$
with $s=\min(s_1,s_2)$. Consequently the class $\fiola = \bigcup _\chi
FIO(\chi ,s)$ is an algebra with respect to the composition of operators. 
\end{theorem}
\begin{proof}
  We write the product $
T^{(1)}T^{(2)}$ as
$$
T^{(1)}T^{(2)} = D_g C_g T^{(1)}T^{(2)}
D_g C_g = D_g (C_g T^{(1)}D_g ) ( C_gT^{(2)} D_g ) C_g \, .
$$
Then $C_g T^{(1)}T^{(2)}
D_g$ is the Gabor matrix of $T^{(1)}T^{(2)}$ with entries
$K_{\mu,\lambda} = $ $\langle T^{(1)}T^{(2)}\pi(\lambda) g,\pi(\mu) g\rangle
$ and $C_g T^{(i)}D_g$, $i=1,2$, is the Gabor matrix of
$T^{(i)}$ with entries $K_{\mu,\lambda}^{(i)}=\langle T^{(i)}\pi(\lambda) g,\pi(\mu)
g\rangle$.  Thus the composition of operators corresponds to the multiplication of
their Gabor matrices.   Using the decay estimates for $K^{(i)}_{\mu, \lambda}$ and
$s=\min (s_1,s_2)$, we estimate the size of the Gabor matrix of
$T^{(1)}T^{(2)}$ as follows:
\begin{align*}
|K_{\mu,\lambda}|& = \sum _{\nu \in \Lambda } K^{(1)} _{\mu ,\nu }
K^{(2)}_{\nu , \lambda } \lesssim\sum_{\nu\in\Lambda} \langle \mu-\chi_1(\nu)\rangle^{-s_1}\langle \nu-\chi_2(\lambda)\rangle^{-s_2}\\
&\lesssim \sum_{\nu\in\Lambda} \langle
\chi_1^{-1}(\mu)-\nu\rangle^{-s} \langle
\nu-\chi_2(\lambda)\rangle^{-s}\, .
\end{align*}
 Since    $v_s \inv = \la
\nu\ra^{-s}\in \ell^1(\Lambda)$ is subconvolutive for $s>2d$, i.e., $v_s\inv \ast
v_s \inv \leq C v_s \inv $~\cite[Lemma 11.1.1 (d)]{book}, the last expression is dominated
by $\langle \chi_1^{-1}(\mu)-\chi_2(\lambda)\rangle^{-s} \asymp
\langle\mu-\chi_1(\chi_2(\lambda))\rangle^{-s}$.

 By
Proposition \ref{P26}, $T^{(1)}$ and $T^{(2)}$ extend to bounded
operators on $L^2(\rd)$, so that  the product $T^{(1)}T^{(2)}$
is well-defined and bounded on $L^2(\rd)$.
\end{proof}

 Finally, we consider the invertibility in the class  $FIO(\chi,s)$
 and show that $\fiola $ is  inverse-closed in $\cB (\lrd )$ for
 $s>2d$.

\begin{theorem} \label{Wineroper} Let $T\in FIO(\chi,s)$ with
  $s>2d$. If  $T$ is invertible on $L^2(\rd)$,
then $T^{-1} \in FIO(\chi^{-1},s)$. Consequently, the algebra $\fiola
$  is inverse-closed in $\cL (\lrd )$.
\end{theorem}
\begin{proof}  We first show that  the adjoint  operator $T^\ast$
  belongs to the class $FIO(\chi \inv, s)$.
  Indeed, since $\chi $ is  bi-Lipschitz, we have
\begin{align*}|\langle T^*\pi(\lambda)g,\pi(\mu) g \rangle|&=|\langle \pi(\lambda)g,T(\pi(\mu) g) \rangle|=|\langle T(\pi(\mu) g,\pi(\lambda)g) \rangle|\\
&\lesssim  \la \lambda-\chi(\mu)\ra^{-s}\asymp\l\
\chi^{-1}(\lambda)-\mu\ra^{-s}.
\end{align*}
Hence, by Theorem \ref{prod1}, the operator $P:=T^\ast T$ is in
$FIO(\mathrm{Id},s)$ and satisfies the estimate $|\langle P\pi
(\lambda )g, \pi (\mu )g\rangle | \lesssim \langle \lambda - \mu
\rangle ^{-s}, \forall \lambda , \mu \in \Lambda $.

 We now  exploit  the
characterization for pseudodifferential operators   contained in
Theorem \ref{CR1} and deduce  that  $P$ is a
pseudodifferential operator with a  symbol in
$S^s_w$. Since $T$ and therefore
$T^\ast$ are invertible on $L^2(\rd)$,  $P$ is also  invertible on
$L^2(\rd)$. 
Now we apply Theorem~\ref{CR2} and conclude that  the inverse
$P^{-1}$ is again a   pseudodifferential
operator with a symbol in $S^s_w$.
Hence $P^{-1}$ is in $FIO(\mathrm{Id},s)$. Finally, using the
algebra property  of Theorem \ref{prod1} once more, we obtain that $T^{-1}=P^{-1}
T^\ast$ is in $FIO(\chi \inv ,s)$ and thus satisfies the estimate
$|\langle T\inv \pi
(\lambda )g, \pi (\mu )g\rangle | \lesssim \langle \chi\inv (\lambda ) - \mu
\rangle ^{-s}, \forall \lambda , \mu \in \Lambda $.
\end{proof}

Combining Theorems~\ref{P26}, \ref{prod1}, and \ref{Wineroper}, we see
that $\fiola $ with $s>2d$  is a Wiener subalgebra of $\cL (\lrd )$
consisting of FIOs.

\section{FIOs of type I}\label{FioI}

The abstract class of Fourier integral operators $FIO(\chi ,s)$ was
defined by decay properties of the Gabor matrix. Our next step is to
relate the Gabor matrix of an operator to the phase and symbol of
concrete FIOs. This step is more technical and resumes our
investigations in ~\cite{CNG,fio1,fio3}.

By the Schwartz' Kernel Theorem every continuous linear
operator $T:\cS(\rd)\to\cS'(\rd)$ can be written as a FIO of type
I with a given phase $\Phi(x,\eta)$  for some symbol
$\sigma(x,\eta)$ in $\cS'(\rdd)$. Hence, if $T$ is a continuous
linear operator $\cS(\rd)\to\cS'(\rd)$ and $\chi$ satisfies {\it
B1}, {\it B2}, \emph{and} {\it B3}, then $T=T_{I,\Phi_\chi,\sigma}$ is a FIO
of type I with symbol $\sigma$ and phase $\Phi_\chi$.
\par

As a first step we formulate the following result.

\begin{proposition}\label{teo1}
Let $T=T_{I,\Phi,\sigma}$ be a FIO of type I  with symbol
$\sigma\in\cS'(\rdd)$ and a phase
 $\Phi$ satisfying A1 and A2. If the Gabor matrix of $T$ satisfies
\begin{equation}\label{uno}
|\langle T \pi(x,\eta) g,\pi(x',\eta')g\rangle|\leq C\langle \nabla_x
\Phi(x',\eta)-\eta',\nabla_\eta \Phi(x',\eta)-x\rangle^{-s} \qquad
x,x',\eta,\eta'\in\R^d \, ,
\end{equation}
for some $C>0$ and  $s\geq 0$, then
$\sigma$ is in the generalized Sj\"ostrand class $ S^s_w(\rdd)$.\par In
particular, if $s>2d$, then  we have $\sigma\in S_w(\rdd)$.
 \end{proposition}
\begin{proof}
The proof uses techniques from \cite{CNG}.
To set up notation,  let  $\Phi _{2,z}$ be  the remainder in the second
 order Taylor expansion  of the phase
 $\Phi $, i.e.,
\begin{equation}\label{phi2}
\Phi_{2,z}(w)=2\sum_{|\alpha|=2}\int_0^1
(1-t)\partial^\alpha\Phi(z+tw)dt \frac{w^\alpha}{\alpha!} \, \qquad
z,w\in\rdd \, ,
\end{equation}
 and set
\begin{equation}\label{psi}
\Psi_z(w)=e^{2\pi i \Phi_{2,z}(w)}
\overline{g}\otimes\widehat{g}(w).
\end{equation}
We recall the fundamental relation between the Gabor matrix of a FIO
and the STFT of its symbol  from \cite[Prop. 3.2]{fio5} and
\cite[Section 6]{fio1}: for $g\in \cS (\rd )$ we have
\[
|\langle T \pi(x,\eta)
g,\pi(x',\eta')g\rangle|=|V_{\Psi_{(x',\eta)}}\sigma
((x',\eta),(\eta'-\nabla_x\Phi(x',\eta),x-\nabla_\eta
\Phi(x',\eta)))| \, .
\]
Writing $u=(x',\eta)$, $v=(\eta',x)$, \eqref{uno}  translates into
\[
|V_{\Psi_u}\sigma(u,v-\nabla\Phi(u))|\leq C\langle v-\nabla \Phi(u)\rangle^{-s},
\]
and then  into the estimate
\begin{equation}\label{due}
\sup_{(u,w)\in\R^{2d}\times \R^{2d}}\langle w\rangle ^s |V_{\Psi_u}\sigma(u,w)|<\infty.
\end{equation}

Now, setting $G=\overline{g}\otimes \widehat{g}\in\cS(\rdd)$, we can write
\begin{align}\label{tre}
V_{G^2}\sigma(u,v)&=\int e^{-2\pi i t v}\sigma(t)\overline{G^2(t-u)}\,dt\nonumber\\
&= \int e^{-2\pi i t v}\sigma(t) e^{-2\pi i \Phi_{2,u}(t-u)}\overline{G(t-u)} e^{2\pi i \Phi_{2,u}(t-u)}\overline{G(t-u)}\, dt\nonumber\\
&= \int e^{-2\pi i t v}\sigma(t) \overline{\Psi_u(t-u)} e^{2\pi i \Phi_{2,u}(t-u)}\overline{G(t-u)}\,dt\nonumber\\
&=\Fur (\sigma T_u\overline{\Psi_u})\ast_v \Fur(T_u(e^{2\pi i
\Phi_{2,u}}\overline{G}))(v)\nonumber\\
&=V_{\Psi_u}\sigma(u,\cdot)\ast \Fur (T_u (e^{2\pi i
\Phi_{2,u}}\overline{G}))(v).
\end{align}
Using \eqref{tre} and the weighted Young inequality $L^\infty_{s}(\rdd)\ast L^1_{{s}}(\rdd)\hookrightarrow L^\infty_s(\rdd)$ we get
\[
\|\sigma\|_{S^s_w}\asymp\sup_u
\| V_{G^2}\sigma(u,\cdot)\|_{L^\infty_s}\lesssim
\sup_u\|V_{\Psi_u}\sigma(u,\cdot)
\|_{L^\infty_s}\sup_u\|\Fur \big(e^{2\pi i
\Phi_{2,u}}\overline{G}\big)\|_{L^1_{s}}.
\]
The first factor in the right-hand side is finite by \eqref{due}. The
second one is finite because the set $\{e^{2\pi i
  \Phi_{2,u}}\overline{G} : u\in \rdd \}$ is bounded in $\cS(\rdd)$, and the embedding
  $\cS \hookrightarrow\Fur L^1_{s}$ is continuous.  This gives $\sigma\in S^s_w$. \par
  The last  statement follows from the inclusion relations for modulation
  spaces in \cite[Proposition 6.5]{F1}, namely,   $M^{\infty,\infty}_{1\otimes
v_s}(\rdd)\hookrightarrow M^{\infty,1}(\rdd)$ if and only if
$s>2d$.
\end{proof}

The next lemma clarifies further the relation between the phase $\Phi
$ and the canonical transformation $\chi$ .

\begin{lemma}\label{tre1}
Consider a phase function
$\Phi$ satisfying {\rm A1},
{\rm A2}, and {\rm A3}.
Then
\begin{equation}\label{3mezzo}
|\nabla_x\Phi(x',\eta)-\eta'|+|\nabla_\eta\Phi(x',\eta)-x|
\asymp |\chi _1(x,\eta)-x'|+| \chi _2(x,\eta)-\eta'| \,  \qquad
\forall x,x',\eta, \eta' \in \rd \, .
\end{equation}
\end{lemma}

\begin{proof}
   The estimate $\gtrsim$ was  already proved in~\cite[Lemma
   3.1]{fio1}).

 For the converse estimate  observe that
   $x=\nabla_\eta\Phi(\chi_1(x,\eta),\eta)$ by definition of $\chi
   _1$, hence 
\begin{equation}\label{4}
|\nabla_\eta \Phi(x',\eta)-x| = |\nabla_\eta \Phi(x',\eta) -
\nabla_\eta \Phi(\chi _1(x,\eta ) ,\eta)| \leq C|x'-\chi_1(x,\eta)|,
\end{equation}
because of assumption (A2) on $\Phi $.
\par Since $\nabla_x \Phi(x',\eta)=\chi_2(\nabla_\eta\Phi(x',\eta),\eta))$, the first term on
the left-hand side of \eqref{3mezzo} can be estimated as
\[
|\nabla_x \Phi(x',\eta)-\eta'|\leq
|\chi_2(\nabla_\eta\Phi(x',\eta),\eta)-\chi_2(x,\eta)|+|\chi_2(x,\eta)-\eta'|
\, .
\]
Finally the  Lipschitz continuity of $\chi $ and \eqref{4} imply  that
\[
|\nabla_x \Phi(x',\eta)-\eta'|\leq|\eta'-\chi_2(x,\eta)|+ C |x'-\chi_1(x,\eta)|.
\]
\end{proof}

The following theorem identifies  the abstract class $FIO(\chi, s)$
with a class of concrete Fourier integral operators and is perhaps the
main result of this paper.

\begin{theorem}\label{caraI}
Fix  $\G(g,\Lambda)$ be a Parseval frame with
$g\in\cS(\rd)$ and let  $s\geq 0$.

Let $T$ be a continuous linear operator $\cS(\rd)\to\cS'(\rd)$ and
$\chi$ a canonical transformation which satisfies {\it B1}, {\it
B2} and {\it B3}.  Then the following properties are
equivalent. \par {\rm (i)} $T=T_{I,\Phi_\chi,\sigma}$ is a FIO of type
I for some $\sigma\in S^s_w$. \par {\rm
(ii)} $F\in FIO(\chi ,s)$. 
\end{theorem}

\begin{proof} The implication $\mathrm{(i)} \Rightarrow \mathrm{(ii)}$  was proved in
  \cite[Theorem 3.3]{CNG}.\par

The implication ${\rm (ii)} \Longrightarrow {\rm (i)}$  follow
immediately  from Proposition~\ref{teo1} and Lemma~\ref{tre1}, since
$\langle (x,\eta ') - \nabla \Phi (x,\eta )\rangle ^{-s} \lesssim
\langle (x',\eta') - \chi (x,\eta )\rangle ^{-s}$. 
\end{proof}

\begin{cor} \label{corhaha}
  Under the same assumptions as in Theorem~\ref{caraI} the following
  statements are equivalent:

{\rm (i)} $T=T_{I,\Phi_\chi,\sigma}$ is a FIO of type
I for some $\sigma\in S^0_{0,0}$.

 \par {\rm
(ii)} $T\in FIO(\chi ,\infty) = \bigcap _{s\geq 0} FIO(\chi ,s)$.
\end{cor}

\begin{remark} \label{tatu} 
\par The corollary should be juxtaposed to Tataru's  characterization
of $FIO(\chi ,\infty )$ in \cite[Theorem 4]{tataru}.  Tataru assumes only conditions  {\it B1}
and {\it B2}, but not {\it B3} on the canonical transformation $\chi
$ and therefore obtains a larger class of FIOs satisfying the decay
condition \eqref{unobis}.  The new insight of Corollary~\ref{corhaha}
is that under the additional assumption {\it B3} \emph{every} FIO
admits the  classical representation \eqref{sei}.
As a byproduct  we see that the integral operator in (14) of \cite{tataru}
possesses  the classical representation, provided that $\chi$ also
satisfies ${\it B3}$.
This observation might be useful for  the
solution of the Cauchy problem in \cite[Section
5]{tataru}, where,  for small time, $\chi$ is a small perturbation
of the identity transformation  and therefore certainly satisfies
{\it B3}.
\end{remark}

Thanks to Theorem \ref{P26}, the FIOs of type I with symbol in
$S^s_w$, with $s>2d$, are bounded on $M^p(\rd)$, $1\leq
p\leq\infty$ (with the obvious modification for $p=\infty$).
Consider now the multiplication operator
\begin{equation}\label{moltiplicatore}T f(x)=e^{\pi i |x|^2}f(x).
\end{equation}
Then $T$ is a FIO of type I having phase $\Phi\phas=|x|^2/2+x
\eta$ and  symbol $\sigma\phas=1$, for every $\phas\in\rdd$.
Observe that $\sigma=1\in S^0_{0,0}$, hence $\sigma\in S_w^s$, for
every $s\geq 0$.
 However, the multiplication operator $T$ is not
bounded on $M^{p,q}$, when $p\not=q$, as proved in
\cite[Proposition 7.1]{fio1}.

We now look at the  Wiener property of the class of FIOs of type I
with symbol in $S^s_w$, with $s>2d$.
As we will see in   Section
\ref{metapoper},   this class  is  not closed under composition, therefore  the Wiener
property must necessarily  involve  FIOs
of type II (see \eqref{sette}).

\begin{theorem} \label{Winerfio}
Let $T$ be a FIO of type I with a tame phase
$\Phi$ and a symbol $\sigma\in S^s_w$, 
with $s>2d$. If $T$ is invertible on $L^2(\rd)$,  then
$T^{-1}$ is a FIO of type II with  same
phase $\Phi$ and  a symbol $\tau\in S^s_w$. 
\end{theorem}
\begin{proof}
Let $\chi $ be the canonical transformation associated to $\Phi
$. Then by  Theorem~\ref{caraI} $T$ belongs to $FIO(\chi ,s)$. As in
the proof of Theorem~\ref{Wineroper} we consider $P=T^*T$ and write
$$
T^{-1}=P^{-1} T^\ast = (T (P^{-1})^\ast)^\ast = (T P\inv )^\ast  \, .
$$
We have already shown that $P$ is in $FIO(\mathrm{Id},s)$ and a \psdo\
with a symbol in $S^s_w$ and that also $P\inv \in FIO(\mathrm{Id},s)$
by the spectral invariance of \psdo s of Theorem~\ref{CR2}.
Now   Theorem~\ref{prod1} implies that $TP\inv \in FIO(\chi, s)$ and
Theorem ~\ref{caraI} implies that $TP\inv $ is a FIO of type I with
tame phase $\Phi $ and a symbol $\rho \in S^s_w$. Since $T\inv =
(TP\inv )^\ast $, $T\inv $ is a FIO of type II with phase $\Phi $ and
the symbol $\tau(x,\eta)=\overline{\rho(\eta,x)}$. An easy
computation as in  \cite[Lemma 2.11]{fio5}) shows that
$\tau\in S^s_w$.
\end{proof}

Although the proof of Theorem~\ref{Winerfio} is short, it combines the main
insights of Sections~3 and~4 and uses the spectral invariance of \psdo
s with symbols in $S^s_w$ (Theorem~\ref{CR2} from \cite{GR}).

\section{Generalized Metaplectic Operators}\label{metapoper}
As an  example, we will  consider  the  class of
$FIO(\chi,s)$ whose phase is a  linear transformation $\chi(z) = \cA
z$ for some invertible matrix $\cA \in GL(2d,\R)$.
Since $\chi $ must preserve the symplectic form (assumption B2), $\cA
$ must be a symplectic matrix, $\cA \in Sp(d,\bR )$. Recall that
the symplectic group is defined  by
$$
\Spnr=\left\{\A\in\Gltwonr:\;^t\!\A J\A=J\right\},
$$
where
$$
J=\begin{pmatrix} 0&-I_d\\I_d&0\end{pmatrix}
\, .
$$

\begin{definition}\label{shsh} Let $\A \in Sp(d,\bR )$ and $s\geq 0$. 
Fix  a Parseval frame  $\G(g,\Lambda)$ with $g\in\cS(\rd)$.  We say that  a continuous linear
operator $T:\,\cS(\rd)\to\cS'(\rd)$ is  a \emph{generalized
  metaplectic operator}, in short, $T\in FIO(\A,s)$,  if its Gabor
matrix satisfies the decay condition
$$
|\langle T \pi(\lambda) g,\pi(\mu)g\rangle|\leq {C}\langle \mu-\cA
\lambda\rangle^{-s},\qquad \forall \lambda,\mu\in \Lambda.
$$
The union $\bigcup _{\cA \in Sp(d,\bR )} FIO(\cA ,s)$ is called the
class of generalized metaplectic operators and denoted by
$FIO(Sp,s)$.
\end{definition}
Since $Sp(d,\bR )$ is a group, 
Theorems \ref{prod1} and \ref{Wineroper} imply the following
statement.
\begin{theorem} For   $s>2d$,  $FIO(Sp,s)$ is a Wiener subalgebra of
$\fiola $.
\end{theorem}

The main examples in $FIO(Sp,s)$ are the operators of the metaplectic
representation of $Sp(d,\bR )$. Given $\cA \in \Spnr $, the
metaplectic operator $\mu (\cA )$ is  defined by the intertwining
relation
\begin{equation}\label{metap}
\pi (\cA z) = c_\cA  \, \mu (\cA ) \pi (z) \mu (\cA )\inv  \quad  \forall
z\in \R ^d \, ,
\end{equation}
where $c_\cA \in \bC , |c_{\cA } | =1$ is a phase factor. The
existence of the metaplectic operators is a consequence of the
Stone-von Neumann theorem  and the irreducibility of the (projective)
representation of $\rdd $ by the \tfs s $\pi (z), z\in \rdd$.
The phase $c_\cA $ can be chosen in such a way that $\mu $
 lifts to a  unitary
representation of the double cover of the symplectic group (which we
will assume in the sequel). For the
group theoretical background and the construction of the  metaplectic
representation we refer to ~\cite{folland89}.
\par
 Let
$\mathcal{A}=\begin{pmatrix} A&B\\C&D\end{pmatrix}\in Sp(d,\R)$ with
$d\times d$ blocks $A,B,C,D$. Then  condition {\it B3}  of Definition \ref{de} is
equivalent to $\det A\not=0$. In this case, $\mu (\cA )$ is explicitly
given by the FIO of type I
\begin{equation}\label{f4}
\mu (\cA )(x)=(\det A)^{-1/2}\int _{\rd }  e^{2\pi i \Phi(x,\eta)}
\hat{f}(\eta)\,d\eta \,
\end{equation}
with  the   phase $\Phi$  
given by
\begin{equation}\label{fase}\Phi(x,\eta)=\frac12  x CA^{-1}x+
\eta  A^{-1} x-\frac12\eta  A^{-1}B\eta \, .
\end{equation}
 (see  \cite[Theorem 4.51]{folland89})



Even without the explicit form of $\mu (\cA )$,  Definition \ref{shsh}
 yields some interesting  information about  the metaplectic
 representation.

\begin{proposition} \label{5.3}
  If $\cA \in \Spnr$, then  $\mu(\A)\in \bigcap _{s\geq 0}  FIO(\A,s)$.
\end{proposition}
 \begin{proof}
 According to Theorem~\ref{cara} it is enough to prove that $\mu(\A)$
 satisfies the continuous  decay condition  \eqref{unobis}.
  Using the definition of $\mu (\cA )$ in
 \eqref{metap} and the covariance property of the STFT~\eqref{eql2}, we
write the Gabor matrix of $\mu (\cA )$ as
\begin{align}
|\langle \mu(\A) \pi(z)g,\pi(w)g\rangle|&=|\langle
 \pi(\A z)\mu(\A)g,\pi(w)g\rangle|\notag\\
 &=|V_{g}(\pi(\A z )\mu(\A)g)(w)|\notag\\
 &=|V_{g}\big(\mu(\A)g\big)\big(w-\A z\big)| \, .
\end{align}
Since both $g\in \cS (\rd )$ and
$\mu(\A)g\in\cS(\rd)$, we also have
 $V_{g}(\mu(\A)g)\in\cS(\rdd)$, e.g., by~\cite[Ch.~11.2.5]{book}.
 This gives
 \begin{equation}
 |\langle \mu(\A) \pi(z)g,\pi(w)g\rangle| \leq C \langle
w-\A z\rangle^{-s} \, , \label{stimam}
 \end{equation}
for every $s\geq0$,  as desired.
\end{proof}

 The following theorem shows that  every
generalized metaplectic operator is a  product of a metaplectic
operator and a classical  pseudodifferential operator.
\begin{theorem}\label{pseudomu}
Let $\cA \in \Spnr $ and $T\in FIO(\A,s)$ with $s>2d$. Then there
exist symbols $\sigma_1, \sigma_2 \in
S^s_w$ with the corresponding pseudodifferential operators   $\sigma_1(x,D)$
and $\sigma_2(x,D)$, such that
\begin{equation}\label{pseudomu1}
T=\sigma_1(x,D)\mu(\A)\quad \mbox{and}\quad
T=\mu(\A)\sigma_2(x,D).
\end{equation}
\end{theorem}
\begin{proof} We prove the factorization $T=\sigma_1(x,D)\mu(\A)$, the other
factorization  is obtained analogously.\par
Since  $\mu(\A)^{-1}=\mu(\A^{-1})$ is in $FIO( \cA \inv ,s)$ by
Proposition~\ref{5.3}, the algebra property of Theorem~\ref{prod1}
implies that
$T\mu(\A^{-1})\in FIO(Id,s)$.   The fundamental
characterization of  pseudodifferential operators of Theorem \ref{CR1}
implies that existence of a symbol
$\sigma _1 \in S^s_w$, such that
$T\mu(\A)^{-1}=\sigma_1(x,D)$, which is what we wanted to show.
\end{proof}

Finally, we  check the counterpart of FIO I and II
for generalized metaplectic operators. \par
 Let
$\mathcal{A}=\begin{pmatrix} A&B\\C&D\end{pmatrix}\in Sp(d,\R)$ with
$\det A\not=0$. As proved in  Theorem~\ref{caraI}, every generalized
metaplectic operator  $T\in FIO(\A,s)$ with
$s>2d$ is a  FIO of type I
\begin{equation}\label{f42}
Tf(x)=(\det \A)^{-1/2}\int e^{2\pi i \Phi(x,\eta)} \sigma\phas
\hat{f}(\eta)\,d\eta,
\end{equation}
with a symbol $\sigma\in S^s_w$ and phase  $\Phi(x,\eta)=\frac12  x CA^{-1}x+
\eta  A^{-1} x-\frac12\eta  A^{-1}B\eta.$

\par We obtain  examples of FIOs of type II by taking adjoints. If
$T\in FIO(\cA ,s)$, then $T^*$ is a FIO of type II.
\par

As is well-known, the generalized
metaplectic operators of type I defined in \eqref{f42} do not enjoy the
algebra property. Consider, for instance, the operators
$T_1=\mu(\A_1)$ and $T_2=\mu(\A_2)$, with
$$
\A_1=\begin{pmatrix} I_d&I_d\\0&I_d\end{pmatrix},\quad
\A_2=\begin{pmatrix} I_d&I_d\\-I_d&0\end{pmatrix}.
$$
Then both  $T_1 $ and $T_2$ are FIOs of type I but their product
$$T_1
T_2=\mu(\A_1)\mu(\A_2)=\mu(\A_1\A_2)=\mu(-J)=\Fur$$ cannot be  a FIO
of type I.
 Indeed, the Fourier transform $\cF = \mu(-J)$ is  an example of a  metaplectic operator that is
 a FIO of  neither type I nor of  type II. Note that in this case assumption
{\it B3} is not satisfied for $\cJ $.

As in Remark~\ref{tatu} we see again that there is a crucial difference between
FIOs satisfying all axioms \emph{B1, B2,} and \emph{ B3} or only
\emph{B1} and \emph{B2}.

\par
\section*{Acknowledgment}
We would like to thank J.~M.~Bony for several discussions on the subject of this paper.

\end{document}